\documentclass[12pt]{article}
\usepackage{amsmath, graphics, latexsym}
\newcommand{\eqdef}{\, =\kern -12.7pt\raise 6pt\hbox{{\tiny\textrm{def}}}\,\,}
\newcommand{\rquote}{\kern-1pt\raisebox{-.4pt}{\scalebox{.7}{\mbox{'}}}}
\newtheorem{theorem}{Theorem}
\newtheorem{cor}{Corollary}
\title{Left-to-right maxima in words and multiset permutations}
\author{Amy N. Myers\\Saint Joseph's University\\Philadelphia, PA 19131
\\{\small \texttt{<amyers@sju.edu>}}
\and
Herbert S. Wilf\\University of Pennsylvania\\Philadelphia, PA 19104
\\{\small \texttt{<wilf@math.upenn.edu>}}}
\begin{document}
\maketitle
\begin{abstract}
We extend classical theorems of R\'enyi by finding the distributions
of the numbers of both weak and strong left-to-right maxima (a.k.a.
\textit{outstanding elements}) in words over a given alphabet and in
permutations of a given multiset.
\end{abstract}

\section{Introduction}
Given a sequence $w=w_1w_2\dots w_n$ of members of a totally ordered
set, we say $j$ is a \emph{strongly outstanding element} of $w$ if
whenever $i < j$ we have $w_i < w_j$. In this case we call $w_j$ a
\emph{strongly outstanding value}.  We say $j$ is a \emph{weakly
outstanding element} if $w_i \leq w_j$ whenever $i < j$, and call
$w_j$ a \emph{weakly outstanding value}. In this paper we will explore
the contexts in which $w$ is a permutation, a multiset permutation, or a
word over some finite alphabet. A famous theorem of R\'{e}nyi
\cite{ren} (see also \cite{lc}) states that the number of
permutations of $[k]$ with $r$ strongly outstanding elements is
equal to the number of such permutations with $r$ cycles, the latter
being given by $\left[\begin{array}{ll}k\\r\end{array}\right]$, the
unsigned Stirling number of the first kind.

In this paper we investigate additional properties of the
outstanding elements and values of permutations, and extend them to
multiset permutations and words on $[k]=\{1,2,\dots,k\}$. An
interesting sidelight to our results is that we obtain a proof of
Gauss's celebrated $_2F_1$ evaluation by comparing two forms of one
of our generating functions, in section \ref{sec:str} below.

\section{Summary of results}
\subsection{Multiset permutations}
For permutations\footnote{Note that for permutations, all outstanding elements and values are
strongly outstanding.} we have the following results.
\begin{theorem}
\label{th:sms} Let $M =\{1^{a_1}, 2^{a_2}, \dots, k^{a_k}\}$ be a
multiset with $N=a_1+a_2+\dots+a_k$, and let $f(M,r)$ denote the
number of permutations of $M$ that contain exactly $r$ strongly
outstanding elements. Then
\begin{equation}
\label{gfso}
F_M(x)\eqdef\sum_r f(M,r)x^r=\frac{(N-1)!a_k x}{a_1!a_2!\dots
a_k!}\prod_{i=1}^{k-1}\left(1+\frac{a_ix}{N-(a_1+a_2+\dots+a_{i})}\right).
\end{equation}
\end{theorem}
\begin{cor}
The generating function for the probability that
a randomly selected permutation of $M$ has exactly $r$ strongly
outstanding elements is
\begin{equation}
\label{pgfso}
P_M(x)=\frac{a_kx}{N}\prod_{i=1}^{k-1}\left(1+\frac{a_ix}{N-(a_1+a_2+\dots+a_{i})}\right),
\end{equation}
and the average number of strongly outstanding elements among
permutations of $M$ is
\begin{equation}
\label{eq:pmpe}
P_M'(1)=\sum_{i=1}^{k}\frac{a_i}{a_i+a_{i+1}+\dots+a_k}.
\end{equation}
\end{cor}

\begin{theorem}
\label{th:woep}
  The generating function for the number of permutations of $M$
that contain exactly $t$ weakly outstanding elements is given by
\begin{equation}
\label{gfwo}
G_M(x)=\phi_{N,a_1}(x)\phi_{N-a_1,a_2}(x)\dots
\phi_{N-a_1-\dots-a_{k-2},a_{k-1}}(x)x^{a_k}
\end{equation}
where
\[\phi_{N,a}(x)=\sum_{m=0}^a{N-m-1\choose a-m}x^m.\]
\end{theorem}
\begin{cor}
The average number of weakly outstanding elements among permutations
of $M$ is
\begin{equation}
\label{eq:avgweak}
\sum_{i=1}^k\frac{a_i}{a_{i+1}+a_{i+2}+\dots+a_k+1}.
\end{equation}
\end{cor}
\begin{cor}
Let $A=\max_i\{a_i\}$.  The amount by which the average number of
weakly outstanding elements exceeds the average number of strongly
outstanding elements is $\le\frac{\pi^2}{6}A(A-1)$.
\end{cor}
\subsection{Words}
The next theorem involves Stirling numbers of the second kind,
denoted by ${n \brace m}$. This is defined as the number of ways to
partition a set of $n$ elements into $m$ nonempty subsets.
\begin{theorem}
\label{th:wordsstrong}
  The number of $n$-letter words over an alphabet of $k$
letters which have exactly $r$ strongly outstanding elements is
given by
\begin{equation}
\label{wordso}
f(n,k,r) = \sum_m{k\choose m}{m\brack r}{n\brace m},
\end{equation}
and the average number of strongly outstanding elements among such
words is
\[H_k=1+\frac12+\frac13+\dots+\frac1k+o(1)\quad(n\to\infty).\]
\end{theorem}
\begin{theorem}
\label{th:wordsweak}
The generating function for the number $g(n,k,t)$ of $n$-letter words
over an alphabet of $k$ letters which have exactly $t$ weakly
outstanding elements is given by
\begin{equation}
\label{wordwo} G_k(n,x)\eqdef \sum_tg(n,k,t)x^t=
\sum_{t=0}^{k-1}(-1)^{k-1-t}(x+t)^n{k-1\choose t}{x+t-1\choose k-1},
\end{equation}
and the average number of weakly
outstanding elements among these words is
\[\frac{n}{k}+H_{k-1}+O\left(\left(\frac{k-1}{k}\right)^n\right).\]
\end{theorem}

Next we introduce the notion of a \emph{template} for words on
$[k]$. A permutation of 5 or more letters \emph{matches the
template} `$YN*YY$', for example, if 1, 4, and 5 are outstanding
elements, 2 is not an outstanding element, and 3 is unconstrained.
For example the permutation 2145763 matches this template. In this
case we think of the letters of the template $Y$, $N$,
* as representing yes, no, and unconstrained, respectively.  We
generalize the $Y$ and $N$ constraints to $S$, $W$, $\overline{S}$,
and $\overline{O}$; where $S$ indicates a strongly outstanding
element, $W$ a weakly outstanding element, $\overline{S}$ indicates
the absence of a strongly outstanding, and $\overline{O}$ the
absence of an outstanding element.  We provide an algorithm for
producing the generating function for the number of words that match
a given template.  When the word is a permutation, we have:
\begin{theorem}
\label{th:tmp1}
Let $\tau$ be a given template of length at most $n$, and let
$\tau_j$ denote the letter that appears in position $j$, counting
from the left, of the template $\tau$.  Since every element of a
permutation is either strongly outstanding or not outstanding, the
letters of $\tau$ are chosen from $\{Y, N, *\}$. The probability
that a permutation of at least $n$ letters matches the template
$\tau$ is
\begin{equation}
\label{notm} \prod_{j:\tau_j={\mathrm{ `N\rquote}}}\left(1-{1\over
{j}}\right)\prod_{j:\tau_j={\mathrm{ `Y\rquote}}}{1\over {j}}
.\end{equation}
\end{theorem}
The corresponding result for words is:\\
\begin{theorem}
\label{th:tmp2}
 Let $\tau$ be a word on $\{S, W, *, \overline{S},
\overline{O}\}$.  Suppose $F(k,\tau,x)=\sum_{k\ge 1}f(k,\tau)x^k$ is
the ordinary generating function for $f(k,\tau)$, the number of
words over the alphabet $[k]$ that match the template $\tau$.
Consider the adjunction of one new symbol, $A \in \{S, W, *,
\overline{S}, \overline{O}\}$, at the right end of $\tau$.  The
generating function, $F(k,\tau A,x)$ can be obtained from
$F(k,\tau,x)$ by applying an operator $\Omega_A$, i.e.,
\[F(k,\tau A,x)=\Omega_AF(k,\tau,x),\]
where
\begin{eqnarray}
\Omega_SF(k,\tau, x) & = & xF(k,\tau, x)/(1-x),\label{omegas}\\
\Omega_WF(k,\tau, x) & = & F(k,\tau, x)/(1-x),\label{omegaw}\\
\Omega_*F(k,\tau, x) & = & x\frac{d}{dx}F(k,\tau, x),\label{omegastar}\\
\Omega_{\overline{S}}F(k,\tau, x) & = &
\left(x\frac{d}{dx}-\frac{x}{1-x}\right)F(k,\tau, x),\ \text{and}\label{omegaovers}\\
\Omega_{\overline{O}}F(k,\tau, x) & = &
\left(x\frac{d}{dx}-\frac{1}{1-x}\right)F(k,\tau,
x)\label{omegaovero}.
\end{eqnarray}
\end{theorem}
\section{Notation}
In the following sections we will count permutations, multiset
permutations, and words that contain a given number of strongly or
weakly outstanding elements.  Let $f(k,r)$ denote the number of
permutations of $[k]=\{1,2,\dots,k\}$ that contain exactly $r$
strongly outstanding elements.  For a multiset
$M=\{1^{a_1},2^{a_2},\dots,k^{a_k}\}$, let $f(M,r)$ denote the
number of permutations of $M$ that contain exactly $r$ strongly
outstanding elements. Finally let $f(n,k,r)$ denote the number of
$n$ letter words on the alphabet $[k]$ that contain exactly $r$
strongly outstanding elements.  When counting weakly outstanding
elements, we use $g$ in place of $f$ and $t$ in place of $r$.

Using the above notation, we define the following generating
functions for strongly outstanding elements: $F_k(x) = \sum_r
f(k,r)x^r$, $F_M(x) = \sum_r f(M,r)x^r$, and $F_k(n,x) = \sum_r
f(n,k,r)x^r$. We define analogous generating functions for weakly
outstanding elements and use $G$ in place of $F$.

\section{Strongly outstanding elements of multiset permutations}
In this section we establish Theorem \ref{th:sms} and its corollaries.

Given a
multiset $M=\{1^{a_1},2^{a_2},\dots ,k^{a_k}\}$, let $N=\sum_ja_j$.
We construct the permutations of $M$ that have exactly $r$ strongly
outstanding elements as follows. We have $N$ slots into which we
will put the $N$ elements of $M$ to make these permutations.

Take the $a_1$ 1's that are available and place them in some
$a_1$-subset of the $N$ slots that are available. There are two
cases now. If the set of slots that we chose for the 1's did not
include the first (leftmost) slot, then we can fill in the remaining
slots with any permutation of the multiset $M/1^{a_1}$ that has
exactly $r$ strongly outstanding elements. On the other hand, if we
did place a 1 into the first slot, then after placing all of the
1's, the remaining slots can be filled in with any permutation of
the multiset $M/1^{a_1}$ that has exactly $r-1$ strongly outstanding
elements.

Recall that $f(M,r)$ denotes the number of permutations of the
multiset $M$ that have exactly $r$ strongly outstanding elements.
The argument in the preceding paragraph shows that
\[f(M,r)=\left({N\choose a_1}-{N-1\choose a_1-1}\right)f(M/1^{a_1},r)+{N-1\choose a_1-1}f(M/1^{a_1},r-1).\]
When we define
\[F_M(x)=\sum_rf(M,r)x^r,\]
we have the recurrence
\begin{eqnarray*}
F_M(x)&=&\left(\left({N\choose a_1}-{N-1\choose
a_1-1}\right)+{N-1\choose a_1-1}x\right)F_{M/1^{a_1}}(x)\\
&=&\left({N-1\choose a_1}+{N-1\choose a_1-1}x\right)F_{M/1^{a_1}}(x)
\end{eqnarray*}
This shows that the generating polynomial resolves into linear
factors over the integers. Indeed we get the explicit form
\begin{eqnarray*}
F_M(x)&=& \left({N-1\choose a_1}+{N-1\choose
a_1-1}x\right)\left({N-a_1-1\choose a_2}+{N-a_1-1\choose
a_2-1}x\right)\dots\\
&=&\prod_{i=1}^k\left({N-a_1-a_2-\dots -a_{i-1}-1\choose
a_i}+{N-a_1-a_2-\dots -a_{i-1}-1\choose a_i-1}x\right)\\
&=&\frac{(N-1)!a_k x}{a_1!a_2!\dots
a_k!}\prod_{i=1}^{k-1}\left(1+\frac{a_ix}{N-(a_1+a_2+\dots+a_{i})}\right)
\end{eqnarray*}
This gives us Theorem \ref{th:sms}. Its corollaries follow by obvious calculations.

Since the generating polynomial has real zeros only, the
probabilities $\{p_r(M)\}$ are unimodal and log concave. By
Darroch's Theorem \cite{dar}, the value of $r$ for which $p_r(M)$ is
maximum differs from $P_M'(1)$ of (\ref{eq:pmpe}) above by at most 1.

\section{Weakly outstanding elements of multiset permutations}
Next we prove Theorem \ref{th:woep}.

Consider $g(M,t)$, the number of permutations of
$M=\{1^{a_1},2^{a_2},\dots, k^{a_k}\}$ that have exactly $t$ weakly
outstanding elements, and $G_M(x)=\sum_tg(M,t)x^t$. To find
$g(M,t)$, suppose the permutation begins with a block of exactly
$m\ge 0$ 1's. Since the value that follows the last 1 in the block
is not available for a 1, there remain $N-m-1$ slots into which the
remaining 1's can be put, in ${N-m-1\choose a_1-m}$ ways. Once all
of the 1's have been placed, if the remaining permutation of the
multiset $M/1^{a_1}$ has exactly $t-m$ weakly outstanding elements
then the whole thing will have $t$ weakly outstanding elements.
Hence we have
\[g(M,t)=\sum_{m\ge 0}{N-m-1\choose a_1-m}g(M/1^{a_1},t-m),\]
if $M/1^{a_1}$ is nonempty, whereas if $M=1^{a_1}$ only, then
$g(M,t)=\delta_{t,a_1}$. If we multiply by $x^t$ and sum on $t$ we
find that
\[G_M(x)=\left\{\sum_{m\ge 0}{N-m-1\choose a_1-m}x^m\right\}G_{M/1^{a_1}}(x),\]
except that if $M=1^{a_1}$ only, then $G_M(x)=x^{a_1}$. So if we
write
\[\phi_{N,a}(x)=\sum_{m=0}^a{N-m-1\choose a-m}x^m,\]
then we have
\begin{equation}
\label{eq:gm}
 G_M(x)=\phi_{N,a_1}(x)\phi_{N-a_1,a_2}(x)\dots
\phi_{N-a_1-\dots-a_{k-2},a_{k-1}}(x)x^{a_k}.
\end{equation}
This gives (\ref{gfwo}).

The two sums $\phi_{N,a}(1)={n\choose a}$ and
$\phi_{N,a}'(1)={n\choose a}/(n-a+1)$ are elementary, and they imply
that $\phi_{N,a}'(1)/\phi_{N,a}(1)=1/(n-a+1)$. Then logarithmic
differentiation of (\ref{eq:gm}) and evaluation at $x=1$ shows that
 the average number of weakly outstanding elements in permutations of
$M$ is
\begin{equation}
\label{eq:wb}
\sum_{i=1}^k\frac{a_i}{a_{i+1}+a_{i+2}+\dots+a_k+1}.
\end{equation}
Let $A=\max_i\{a_i\}$. If we compare (\ref{eq:wb}) and (\ref{eq:pmpe})
we find that the amount by which the average number of weakly
outstanding elements exceeds the average number of strongly
outstanding elements is
\begin{eqnarray*}
\sum_{i=1}^k&&\kern-25pt\left(\frac{a_i}{a_{i+1}+a_{i+2}+\dots+a_k+1}
-\frac{a_i}{a_i+a_{i+1}+\dots+a_k}\right)\\
\qquad&=&\sum_{i=1}^k\frac{a_i(a_i-1)}{(a_{i+1}+a_{i+2}+\dots+a_k+1)(a_i+a_{i+1}+\dots+a_k)}\\
\qquad&\le& A(A-1)\sum_{i=1}^k\frac{1}{(k-i+1)^2}\le \frac{\pi^2}{6}A(A-1).
\end{eqnarray*}
This estimate is best possible when $A=1$, i.e., when every element
occurs just once.

\section{Strongly outstanding elements of words}
\label{sec:str}
 Next we investigate $f(n,k,r)$, the number of
$n$-letter words over an alphabet of $k$ letters that have exactly
$r$ strongly outstanding elements.  In so doing we will prove
Theorem \ref{th:wordsstrong}.

Note first that the number of strongly outstanding elements of such
a word depends only on the permutation of the distinct letters
appearing in the word that is achieved by the \textit{first}
appearances of each of those letters, because a value $j$ can be
strongly outstanding in a word $w$ only if it is the first (i.e.,
leftmost) occurrence of $j$ in $w$.

Hence associated with each $n$-letter word $w$ over an alphabet of
$k$ letters which has exactly $r$ strongly outstanding elements
there is a triple $(S,{\cal P},\sigma)$ consisting of
\begin{enumerate}
\item a subset $S\subseteq [k]$, which is the set of all of the distinct
 letters that actually appear in $w$, and
\item a partition ${\cal P}$ of the set $[n]$ into $m=|S|$ classes,
namely the $i$th class of ${\cal P}$ consists
 of the set of positions in the word $w$ that contain the $i$th letter
  of $S$, and
\item a permutation $\sigma\in S_m$, $m=|S|$, which is the sequence of \textit{first} appearances in $w$ of each of the $k$ letters that occur in $w$. $\sigma$ will have $r$ strongly outstanding elements.
\end{enumerate}

Conversely, if we are given such a triple $(S,{\cal P},\sigma)$, we
uniquely construct an $n$-letter word $w$ over $[k]$ with exactly
$r$ strongly outstanding elements as follows.

First arrange the classes of the partition ${\cal P}$ in ascending
order of their smallest elements. Then permute the set $S$ according
to the permutation $\sigma$, yielding a list $\tilde{S}$. In all of
the positions of $w$ that are described by the first class of the
partition ${\cal P}$ (i.e., the class in which the letter `1' lives)
we put the first letter of $\tilde{S}$, etc., to obtain the required
word $w$.

Thus the number of words that we are counting is equal to the number
of these triples, viz.
\begin{equation}
\label{eq:fnkr}
 f(n,k,r)=\sum_m{k\choose m}{m\brack r}{n\brace m}.
\end{equation}
It is noteworthy that three flavors of ``Pascal-triangle-like''
numbers occur in this formula.

Let $\rho(n,k)$ denote the average number of strongly outstanding
elements among the $n$-letter words that can be formed from an
alphabet of $k$ letters. To find $\rho(n,k)$ for large $n$, we have
first that \begin{eqnarray*} \sum_rf(n,k,r)&=&\sum_m{k\choose
m}{n\brace m}\sum_r{m\brack r}= \sum_m{k\choose m}{n\brace
m}m!\\
&\sim& \sum_m{k\choose m}{n\brace m}m!\sim \sum_m{k\choose m}n^m\sim
n^k,
\end{eqnarray*}
for large $n$, where we have used the facts that ${n\brace m}\sim
n^m/m!$ and $\sum_r{m\brack r}=m!$. Similarly,
\begin{eqnarray*}\sum_rrf(n,k,r)&=&\sum_m{k\choose m}{n\brace m}\sum_rr{m\brack
r}= \sum_m{k\choose m}{n\brace m}m!H_m\\
&\sim& \sum_m{k\choose m}{n\brace m}m!H_m\sim \sum_m{k\choose
m}n^mH_m\sim n^kH_k,\end{eqnarray*}
 where we have used the additional fact that
$\sum_rr{m\brack r}=m!H_m$. If we divide these last two equations we
find that
\[\lim_{n\to\infty}\rho(n,k)=H_k=1+\frac12+\frac13+\dots+\frac1k.\]
The proof of Theorem \ref{th:wordsstrong} is complete. $\Box$

From eq. (\ref{eq:fnkr}) we can use the standard generating
functions for the two kinds of Stirling numbers to show that
\[\sum_{r,n}f(n,k,r)y^rt^{n-r}=\prod_{j=1}^k\left(1+\frac{y}{1-jt}\right).\]
But from (\ref{eq:fnkr}) we also have
\begin{eqnarray*}
\sum_{r,n}f(n,k,r)y^rt^{n-r}&=&\sum_{\ell,r}{k\choose \ell}{\ell\brack r}\left(\frac{y}{t}\right)^r\sum_n{n\brace \ell}t^n\\
&=&\sum_{\ell,r}{k\choose \ell}{\ell\brack r}\left(\frac{y}{t}\right)^r\frac{t^{\ell}}{(1-t)(1-2t)\dots (1-\ell t)}\\
&=&\sum_{\ell}{k\choose \ell}\prod_{j=0}^{\ell-1}\left(\frac{y}{t}+j\right)\frac{t^{\ell}}{(1-t)(1-2t)\dots (1-\ell t)}\\
&=&\sum_{\ell}{k\choose \ell}\prod_{j=1}^{\ell}\frac{y+(j-1)t}{1-jt}.
\end{eqnarray*}

Comparison of these two evaluations shows that we have found and
proved the following identity:
\begin{equation}
\label{eq:ident}
\sum_{\ell}{k\choose \ell}\prod_{j=1}^{\ell}\frac{y+(j-1)t}{1-jt}=\prod_{j=1}^k\left(1+\frac{y}{1-jt}\right).
\end{equation}
But Gauss had done it earlier, since it is the evaluation of his
well known
\[ _2F_1\left[ \begin{array}{cc}
    -k & y/t \\
    1-1/t&
     \end{array}\bigg|\,1 \right].\]

\section{A calculus of templates}
\subsection{Templates and permutations}
In this section we prove Theorems \ref{th:tmp1} and \ref{th:tmp2}.
  We begin by establishing
equation (\ref{notm}), which appeared first in \cite{wil}.

This result is a generalization of a theorem of R. V. Kadison
\cite{rvk}, who discovered the case where the template is `NN$\cdots
$NY', and proved it by the sieve method. To make this paper
self-contained, we include a proof of (\ref{notm}).

Our proof is by induction on $n$. Suppose it has been proved that,
for all templates $\tau$ of at most $n-1$ letters, the number of
permutations of $n-1$ letters that match $\tau$ is correctly given
by (\ref{notm}), and let $\tau$ be some template of $\le n$ letters.
If in fact the length of $\tau$ is $<n$ then the formula
(\ref{notm}) gives the same result as it did when applied to
$(n-1)$-permutations, which is the correct probability.

In the case where the length of $\tau$ is $n$ and the rightmost
letter of $\tau$ is `Y', every matching permutation $\sigma$ must
have $\sigma (n)=n$. Hence the number of matching permutations is
$(n-1)!p_{n-1}(\tau')$, where $\tau'$ consists of the first $n-1$
letters of $\tau$, which is equal to $n!p_n(\tau)$, proving the
result in this case.

In the last case, where the length of $\tau$ is $n$ and the
rightmost letter of $\tau$ is `N', the probability of a permutation
match must be $p_{n-1}(\tau')(1-1/n)$, since this case and the
preceding one are exhaustive of the possibilities and the preceding
one had a probability of $p_{n-1}(\tau')/n$. But this agrees with
the formula (\ref{notm}) for this case, completing the proof of the
theorem.

\subsection{Templates on words}
Next we include results for words that are analogous to those we
found for permutations.  A preview of the kind of results that we
will get is the following. Suppose $F(k,\tau,x)=\sum_{k\geq 1}
f(k,\tau)x^k,$ where $f(k,\tau)$ denotes the number of words over
the alphabet $[k]$ that match the template $\tau$. Then consider the
adjunction of one new symbol, let's call it $A$, at the right end of
$\tau$. Then we will show that the new generating function,
$F(k,\tau A,x)$ can be obtained from $F(k,\tau,x)$ by applying a
certain operator $\Omega_A$. That is,
\[F(k,\tau A,x)=\Omega_AF(k,\tau,x).\]
The operator $\Omega_A$ will depend only on the letter $A$ that is
being adjoined to the template $\tau$.

Hence, to find the generating function for a complete template
$\tau$, we begin with $F(k,\emptyset,x)=1/(1-x)$, and we read the
template $\tau$ from left to right. Corresponding to each letter in
$\tau$ we apply the appropriate operator $\Omega$. When we have
finished scanning the entire template the result will be the desired
generating function for $\tau$.

The letters that we will allow in a template
$\tau=\tau_1\tau_2\dots\tau_l$ are
$\{S,W,O,*,\overline{S},\overline{O}\}$. Their meanings are that if
$w=w_1\dots w_l$ is a word of length $l$ over the alphabet $[k]$
then for $w$ to match the template it must be that $w_i$ is
\begin{enumerate}
\item a strongly outstanding value whenever $\tau_i=S$, or
\item a weakly outstanding value whenever $\tau_i=W$, or
\item unrestricted whenever $\tau_i=*$, or
\item not a strongly outstanding value whenever $\tau_i=\overline{S}$, or
\item neither a weakly nor a strongly outstanding value whenever $\tau_i=\overline{O}$.
\end{enumerate}
Let's consider what happens to the count of matching words when we
adjoin one of these letters to a template whose counting function is
known. Let $f(k,\tau)$ denote the number of words of length $l=$
length$(\tau)$, on the alphabet $[k]$ that match the template
$\tau$.
\begin{enumerate}
\item If $w$ is one of the words counted by $f(k,\tau S)$, and if we delete its last letter, we obtain one of the words that is counted by $f(i,\tau)$ for some $1\le i<k$, and consequently
\[f(k,\tau S)=\sum_{i=1}^{k-1}f(i,\tau).\]
If $F(k,\tau,x)=\sum_{k\ge 1}f(k,\tau)x^k$, then we have
\[F(k,\tau S,x)=\frac{x}{1-x}F(k,\tau ,x).\]
Thus we have found the operator $\Omega_S$, and it is defined by
$\Omega_SF(x)=xF(x)/(1-x)$.  This is equation (\ref{omegas}).
\item Similarly, if $w$ is one of the words counted by $f(k,\tau W)$, and if we delete its last letter, we obtain one of the words that is counted by $f(i,\tau)$ for some $1\le i\le k$, and consequently
\[f(k,\tau W)=\sum_{i=1}^{k}f(i,\tau).\]
If $F(k,\tau,x)=\sum_{k\ge 1}f(k,\tau)x^k$, then we have
\[F(k,\tau W,x)=\frac{F(k,\tau ,x)}{1-x}.\]
Thus we have found the operator $\Omega_W$, and it is defined by
$\Omega_WF(x)=F(x)/(1-x)$.  This is (\ref{omegaw}).
\end{enumerate}
Since the argument in each case is easy and similar to the above we
will simply list the remaining three operators, equations
(\ref{omegastar}), (\ref{omegaovers}), and (\ref{omegaovero}), as
follows:
\begin{eqnarray*}
\Omega_*F(x)&=&x\frac{d}{dx}F(x)\\
\Omega_{\overline{S}}F(x)&=&\left(x\frac{d}{dx}-\frac{x}{1-x}\right)F(x)\\
\Omega_{\overline{O}}F(x)&=&\left(x\frac{d}{dx}-\frac{1}{1-x}\right)F(x)
\end{eqnarray*}

The successive applications of these operators can be started with
the generating function for the empty template,
\[F(\emptyset,x)=\frac{1}{1-x}.\]

Thus to find the generating function for some given template $\tau$,
we begin with the function $1/(1-x)$, and then we read one letter at
a time from $\tau$, from left to right, and apply the appropriate
one of the five operators that are defined above.

As an example, how many words of length 3 over the alphabet $[k]$
match the template $\tau=S*\overline{S}$? This is the coefficient of
$x^k$ in
\begin{eqnarray*}
F(S*\overline{S},x)&=&\Omega_{\overline{S}}\Omega_*\Omega_S\frac{1}{1-x}\\
&=&\left(x\frac{d}{dx}-\frac{x}{1-x}\right)\left(x\frac{d}{dx}\right)
\left(\frac{x}{1-x}\right)\frac{1}{1-x}\\
&=&\frac{x+3x^2}{(1-x)^4}
\end{eqnarray*}
The required number of words over a $k$ letter alphabet that match
the template $\tau$ is the coefficient of $x^k$ in the above. Since
these generating functions will always be of the form
$P(t)/(1-x)^r$, with $P$ a polynomial, we note for ready reference
that
\[[x^k]\left\{\frac{\sum_ja_jx^j}{(1-x)^r}\right\}=\sum_ja_j{r+k-j-1\choose
r- 1}.\] In the example above we have $r=4$, $a_1=1$, $a_2=3$, so
\[f(k,S*\overline{S})={k+2\choose 3}+3{k+1\choose 3}.\]

\section{Weakly outstanding elements of words}
Finally we prove Theorem \ref{th:wordsweak}.

 Recall that $g(n,k,t)$ is the number
of $n$-letter words over the alphabet $[k]$ which have exactly $t$
weakly outstanding elements. Consider just those words $w$ that
contain exactly $m$ 1's, $m<n$. If $w$ begins with a block of
exactly $l$ 1's, $0\le l\le m$ then by deleting all $m$ of the 1's
in $w$ we find that the remaining word is one with $n-m$ letters
over an alphabet of $k-1$ letters and it has exactly $t-l$ weakly
outstanding elements. Thus we have the recurrence
\[
g(n,k,t)=\sum_{m=0}^{n-1}{n-l-1\choose
m-l}g(n-m,k-1,t-l)+\delta_{t,n}.
\]
When we set $G_k(n,x)=\sum_tg(n,k,t)x^t$, we find that
\begin{equation}
\label{eq:f2rec} G_k(n,x)=\sum_{m=0}^{n-1}\sum_{l=0}^k{n-l-1\choose
m-l}x^lG_{k-1}(n-m,x)+x^n.\qquad(n,k\ge 1;\, G_0(n,x)=0)
\end{equation}
To discuss this recurrence, let $\Delta=\Delta_x$ be the usual
forward difference operator on $x$, i.e., $\Delta_x
f(x)=f(x+1)-f(x)$. Then from the recurrence above we discover that
\[G_1(n,x)=x^n;\ G_2(n,x)= \Delta_x x^n(x-1);\ G_3(n,x)=\Delta_x^2x^n{x-1\choose 2}.\]
This leads to the conjecture that
\begin{eqnarray*}
G_k(n,x)&=&\Delta_x^{k-1}\left\{x^n{x-1\choose k-1}\right\}\\
&=&\sum_{t=0}^{k-1}(-1)^{k-1-t}(x+t)^n{k-1\choose t}{x+t-1\choose
k-1}.
\end{eqnarray*}
To prove this it would suffice to show that the function $G_k(n,x)$
above satisfies the recurrence (\ref{eq:f2rec}). If we substitute
the conjectured form of $G_k$ into the right side of
(\ref{eq:f2rec}) we find that the sums over $l$ and $m$ can easily
be done, and the identity to be proved now reads as
\begin{eqnarray*}
&&\sum_{t=0}^{k-2}(-1)^{k-t}{k-2\choose t}{x+t-1\choose k-2}\frac{(x+t)}{t+1}((x+t+1)^n-x^n)+x^n\\
&&\qquad=\sum_{t=0}^{k-1}(-1)^{k-1-t}(x+t)^n{k-1\choose
t}{x+t-1\choose k-1}.
\end{eqnarray*}
If we replace the dummy index of summation $t$ by $t-1$ on the left
side, the identity to be proved becomes
\begin{eqnarray*}
&&\sum_{t=1}^{k-1}(-1)^{k-t-1}{k-2\choose t-1}{x+t-2\choose k-2}\frac{(x+t-1)}{t}((x+t)^n-x^n)+x^n\\
&&\qquad=\sum_{t=0}^{k-1}(-1)^{k-1-t}(x+t)^n{k-1\choose
t}{x+t-1\choose k-1}.
\end{eqnarray*}
It is now trivial to check that for $1\le t\le k-1$, the coefficient
of $(x+t)^n$ on the left side is equal to that coefficient on the
right. If we cancel those terms and divide out a factor $x^n$ from
what remains, the identity to be proved becomes
\[-\sum_{t=1}^{k-1}(-1)^{k-t-1}{k-2\choose t-1}{x+t-2\choose k-2}\frac{(x+t-1)}{t}+1=(-1)^{k-1}{x-1\choose k-1}.\]
The sum that appears above is a special case of Gauss's original
$_2F_1[a,b;c|1]$ evaluation, and the proof of (\ref{wordwo}) is
complete.  A straightforward calculation now shows that the average
number of weakly outstanding elements among $n$-letter words over
the alphabet $[k]$ is
\[=\frac{n}{k}+H_{k-1}+O\left(\left(\frac{k-1}{k}\right)^n\right).\]

\section{Related results}
In the literature, the outstanding elements and values of sequences
go by various names: \emph{\'{e}l\'{e}ments salients} (R\'{e}nyi),
\emph{left-to-right maxima}, \emph{records}, and others; and appear
in several different contexts and results, for example:

It is well known that the probability of obtaining at least $r$
strongly outstanding elements in a sequence $X_1,X_2,\dots ,X_n$ of
$n$ independent, identically distributed, continuous random
variables approaches 1 as $n\rightarrow\infty$.  (Glick's survey
\cite{gl}, for example, contains this result and those that follow
in this paragraph and the next.) Let $Y_i = 1$ if $i$ is a strongly
outstanding element of such a sequence, and 0 otherwise. The
expected value is $E[Y_i]=1/i$ and the variance is $V[Y_i]=1/i -
1/i^2$. The number of strongly outstanding elements in a sequence of
continuous iid random variables is therefore $R_i = \sum_i Y_i$ with
expectation $E[R_i]=\sum_{i=1}^n 1/i$ and variance
$V[R_i]=\sum_{i=1}^n 1/i - \sum_{i=1}^n 1/i^2$. Note $\sum_{i=1}^n
1/i - \ln(n) \rightarrow$ Euler's constant $=.5772\dots$, and
$\sum_{i=1}^n 1/i^2 \rightarrow \pi^2/6 = 1.6449\dots$ for
$n\rightarrow \infty.$

Let $N_r$ denote the $r^{\text{th}}$ outstanding element in a
sequence of $n$ iid continuous random variables.  Then $N_1 = 1$,
$E[N_r] = \infty$ for $r\geq 2$, and $E[N_{r+1} - N_r] = \infty$ as
well.  The probability $P[N_2 = i_2, N_3 = i_3, \dots , N_r = i_r] =
1/(i_2 -1)(i_3-1)\dots (i_r-1)i_r$ for $1 < i_2 < i_3 < \dots <i_r$.
The probability that an iid sequence of $n$ continuous random
variables has exactly $r$ strongly outstanding elements is
$\left[\begin{array}{ll}k\\r\end{array}\right]/(r-1)!$ for large
$n$.

Chern and Hwang \cite{ch} consider the number $f_{n,k}$ of $k$
consecutive records (strongly outstanding elements) in a sequence of
$n$ iid continuous random variables.  They improve upon known
results for the limiting distribution of $f_{n,2}$.  In particular
they show $f_{n,k}$ is asymptotically Poisson for $k=1,2$, and this
is not the case for $k\geq 3$.  They give the probability generating
function for $f_{n,2}$, and observe that the distribution of
$f_{n,2}$ is identical to that for the number of fixed points $j$ in
a random permutation of $[n]$ for $1 \leq j < n$.   They give a
recurrence for the probability generating function for $f_{n,k}$,
and compute the mean and variance for this number.

Are there similar results for discrete distributions?  Prodinger
\cite{prod} considers left-to-right maxima in both the strict
(strongly outstanding) and loose (weakly outstanding) senses for
geometric random variables. He finds the generating function for the
probability that a sequence of $n$ independent geometric random
variables (each with probability $pq^{v-1}$ of taking the positive
integer value $v$, where $q = p-1$) has $k$ strict left-to-right
maxima. This probability is the coefficient of $z^ny^k$ in
$$
F(z,y) = \prod_{k \geq 1}
\left(1+\frac{yzpq^{k-1}}{1-z(1-q^k)}\right).
$$
For loose left-to-right maxima the analogous generating function is
$$
\prod_{k\geq 0}\frac{1-z(1-q^k)}{1-z+zq^k(1-py)}.
$$
In this paper, Prodinger also finds the asymptotic expansions for
both the expected numbers of strict and loose left-to-right maxima,
and the variances for these numbers.  He does all of the above for
uniform random variables as well.

Knopfmacher and Prodinger \cite{kp} consider the \emph{value} and
\emph{position} of the $r^{\text{th}}$ left-to-right maximum for $n$
geometric random variables. The \emph{position} is the
$r^{\text{th}}$ strongly outstanding element, and the \emph{value}
is that taken by the random variable in that position (again the
value $v$ is taken with probability $pq^{v-1}$, where $q=1-p$). For
the $r^{\text{th}}$ strong left-to-right maximum, the asymptotic
formulas for value and position are $\frac{r}{p}$ and
$\frac{1}{(r-1)!}\left(\frac{p}{q}\log_{\frac{1}{q}}n\right)^{r-1}$,
respectively. For the \emph{r}th weak left-to-right maximum, the
value and position are asymptotically $\frac{rq}{p}$ and
$\frac{1}{(r-1)!}\left(p\log_{\frac{1}{q}}n\right)^{r-1}$.  These
results are obtained by first computing the relevant generating
functions.

A number of additional properties of outstanding elements of permutations
 are in Wilf \cite{wil}.

Key \cite{key} describes the asymptotic behavior of the number of
records (strongly outstanding elements) and weak records (weakly
outstanding elements) that occur in an iid sequence of integer
valued random variables.

When the sequence $w_1, w_2, \dots, w_n$ under consideration is a
random permutation on $n$ letters, then the expected number of
strongly outstanding elments is the \emph{n}th harmonic number $H_n
= 1 + \frac{1}{2} + \frac{1}{2} + \dots + \frac{1}{n}$. The variance
is $H_n = H_n^{(2)}$, where
$H_n^{(2)}=1+\frac{1}{4}+\dots+\frac{1}{n^2}$ denotes the \emph{n}th
harmonic number of the second order.

Banderier, Mehlhorn, and Brier \cite{bmb} show that the average
number of left-to-right maxima (strongly outstanding elements) in a
\emph{partial permutation} is $\log(pn) + \gamma +
2\frac{1-p}{p}+\left(\frac{1}{2}+\frac{2(1-p)}{p^2}\right)\frac{1}{n}+O\left(\frac{1}{n^2}\right)$
where $\gamma = .5772\dots$ is Euler's constant.  To obtain a
\emph{partial permutation}, we begin with the sequence $1,2,\dots
,k$ and select each element with probability $p$; then take one of
the $(pn)!$ permutations of $[pn]$ uniformily at random and let it
act on the selected elements, while the nonselected elements stay in
place.

Foata and Han \cite{fh} investigate the (right-to-left) \emph{lower
records} of \emph{signed permutations}.  A \emph{signed permutation}
is a word $w=w_1w_2\dots w_n$ for which the letters $w_i$ are
positive or negative integers, and $|w_1||w_2|\dots |w_n|$ is a
permutation of $\{1,2,\dots,n\}$.  A \emph{lower record} of such a
word is a letter $w_i$ such that $w_i < w_j$ for all $j$ with $i+1
\leq j \leq n$. The authors consider the \emph{signed subword}
obtained by reading the lower records of $w$ from left to right; and
derive generating functions for signed permutations in terms of
signed subwords, numbers of positive and negative letters in signed
subwords, and other statistics.


\newpage


\begin{thebibliography}{hhh}
\bibitem{lc} Louis Comtet, Advanced Combinatorics, Reidel, 1974.
\bibitem{bmb} Cyril Banderier, Kurt Mehlhorn, and Rene Beier,
Smoothed analysis of three combinatorial algorithms, Mathematical
Foundations of Computer Science 2003: 28th International Symposium,
Proceedings, Lecture Notes in Computer Science 2747, Springer, 2003,
198--207.
\bibitem{ch} Hua-Huai Chern and Hsien-Kuei Hwang, Limit distribution
of the number of consecutive records, Random Structures Algorithms
\textbf{26} no.4 (2005), 404--417.
\bibitem{dar} J. N. Darroch, On the distribution of the number of successes in
independent trials, Ann. Math. Stat., \textbf{35} (1964),
1317--1321.
\bibitem{fh} Dominique Foata and Guo-Niu Han, Signed words and permutations II: The Euler-Mahonian polynomials,
Electron. J. Combin.  \textbf{11} no. 2 (2004/06), Research Paper
22, 18 pp. (electronic).
\bibitem{gl} Ned Glick, Breaking records and breaking boards,
American Mathematical Monthly, \textbf{85} (1978), 2--26.
\bibitem{key} Eric S. Key, On the number of recordes in an iid
discrete sequence, Journal of Theoretical Probability, \textbf{18}
no. 1 (2005), 99--108.
\bibitem{rvk} Richard V. Kadison, Strategies in the secretary problem, Expositiones Mathematicae {\bf 12} (1994), 125-144.
\bibitem{kp} Arnold Knopfmacher and Helmut Prodinger, Combinatorics of geometrically distributed random variables:
Value and position of the $r^{th}$ left-to-right maximum.  Discrete
Math.  \textbf{226} no. 1-3 (2001), 255--267.
\bibitem{prod} Helmut Prodinger, Combinatorics of geometrically
distributed random variables:  Left-to-right maxima.  Discrete
Mathematics, \textbf{153} (1996), 253--270.
\bibitem{ren} A. R\' enyi, Th\' eorie des \' el\' ements saillants d'une suite d'observations, in {\it Colloquium Aarhus}, 1962, 104-117.
\bibitem{wil} Herbert S. Wilf, On the outstanding elements of permutations, \texttt{<http://www.cis.upenn.edu/}$~$\texttt{wilf>}, 1995.
\end{thebibliography}
\end{document}